\documentclass[1p]{elsarticle}

\usepackage{lineno,hyperref}
\usepackage{amsmath, amssymb}
\newcommand{\sq}{\qquad $\blacksquare$}
\modulolinenumbers[5]

\makeatletter
\def\ps@pprintTitle{%
 \let\@oddhead\@empty
 \let\@evenhead\@empty
 \def\@oddfoot{}%
 \let\@evenfoot\@oddfoot}
\makeatother

\newtheorem{thm}{Theorem}
\newdefinition{defn}[thm]{Definition}

\begin{document}

\begin{frontmatter}

\title{Hamilton Cycles in Double Generalized Petersen Graphs}

\author{Yutaro Sakamoto\fnref{myfootnote}}
\address{Department of Informatics, The University of Electro-Communications, 1-5-1 Chofugaoka, Chofu, Tokyo 182-8585, Japan}
\ead{y-sakamoto@uec.ac.jp}




\begin{abstract}
  Watkins (1969) first introduced the generalized Petersen graphs (GPGs) by modifying Petersen graph.
  Zhou and Feng (2012) modified GPGs and introduced the double generalized Petersen graphs (DGPGs).
  Kutnar and Petecki (2016) proved that DGPGs are Hamiltonian in special cases
  and conjectured that all DGPGs are Hamiltonian.
  In this paper, we construct Hamilton cycles in all DGPGs.
\end{abstract}


\end{frontmatter}

\linenumbers

\section{Introduction}

\par
In \cite{Wat} Watkins (1969) first introduced the GPGs to discover Tait coloring of the graphs
and Castagna and Prins (1972) proved Watkins' conjecture about GPGs in \cite{Cas}.
After they introduced GPGs, some properties of GPGs have been studied.
For instance, Alspach (1983) determined which GPGs have a Hamilton cycle in \cite{Als}.
Fu, Yang and Jiang (2009) studied the domination number of GPGs in \cite{Fu}.

\par
Now we define double generalized Petersen graphs DP$(n, t)$ (DGPG for short) as follows.
\begin{defn}
  Let $n$ and $t$ be integers that satisfy $n \geq 3$ and $2 \leq 2t < n$.
  The double generalized Petersen graph DP$(n, t)$ is an undirected simple graph with vertex set $V$ and edge set $E$, where
  \begin{align*}
    V &= \{x_i, u_i, v_i, y_i \mid i \in \mathbb{Z}_n \}, \\
    E &= \{x_ix_{i+1}, y_iy_{i+1}, x_iu_i, y_iv_i, u_iv_{i+t}, v_iu_{i+t} \mid i \in \mathbb{Z}_n \}
  \end{align*}
\end{defn}
Note that $\mathbb{Z}_n$ denotes a set of integers $\mathbb{Z}/n\mathbb{Z}$ throughout this paper.

\begin{figure}[!tpb]
  \begin{center}
    \includegraphics[height=0.4\textheight]{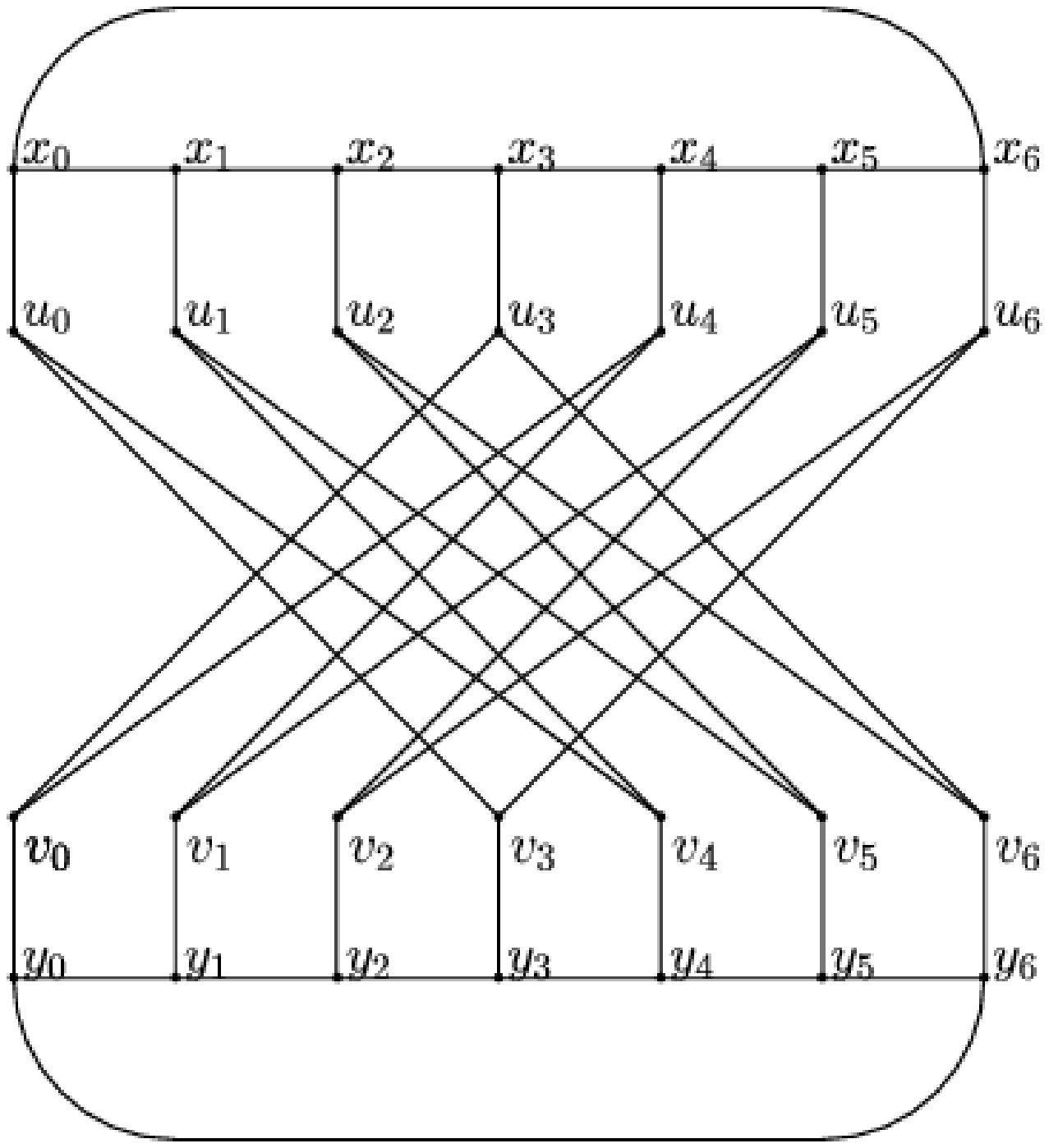}
    \caption{DP$(7, 3)$}
    \label{fig:dp73}
  \end{center}
\end{figure}

Zhou and Feng (2012) first introduced the double generalized Petersen graphs by modifying the generalized Petersen graphs in \cite{Jin1}.
In \cite{Jin2}, Zhou and Feng (2014) determined all non-Cayley vertex-transitive graphs and all vertex-transitive graphs among DGPGs.
From their result, Kutnar and Petecki (2016) gave the complete classification of automorphism groups of DGPGs in \cite{Kutnar}.
They also proved that DP$(n, t)$ is Hamiltonian if $n$ is even or $n$ is odd and the greatest common divisor of $n$ and $t$ equals to $1$ in \cite{Kutnar}.
In addition, a computer-assisted search verified that DP$(n,t)$ have Hamilton cycles for all $n \leq 31$ 
and they conjectured that all DP$(n, t)$ are Hamiltonian.
This paper gives the following theorem.

\begin{thm}\label{thm}
  All DP$(n, t)$ are Hamiltonian.
\end{thm}

\section{Preliminaries}
  As mentioned in the previous section, $\mathbb{Z}_n$ denotes a set of integers $\mathbb{Z}/n\mathbb{Z}$.
  A sequence of vertices $w_0 w_1 w_2 \dots w_n$ denotes a path in a graph.
  A path whose end points are the same vertex is called a cycle.
  $V(G)$ denotes the vertex set of a graph $G$.
  Let $G$ be an arbitrary subgraph of DP$(n, t)$.
  We define functions $V_x, V_y, V_u, V_v$ as follows.
  \begin{align*}
    V_x(G) &= V(G) \cap \{x_i \mid i \in \mathbb{Z}_n\}\\
    V_y(G) &= V(G) \cap \{y_i \mid i \in \mathbb{Z}_n\}\\
    V_u(G) &= V(G) \cap \{u_i \mid i \in \mathbb{Z}_n\}\\
    V_v(G) &= V(G) \cap \{v_i \mid i \in \mathbb{Z}_n\}\\
  \end{align*}

\section{The construction of Hamilton cycles in DP$(n, t)$}
  We assume that $n$ is even. In this case, Kutnar and Petecki showed that all DP$(n,t)$ are Hamiltonian in \cite{Kutnar}.
  Observe that there exist paths $X_i$  for all $i \in \mathbb{Z}_{n/2}$.
  \begin{align*}
    X_i : u_{2i} x_{2i} x_{2i+1} u_{2i+1} v_{2i+1-t} y_{2i+1-t} y_{2i+2-t} v_{2i+2-t} u_{2(i+1)}
  \end{align*}
  Joining all of the paths gives a Hamilton cycle in DP$(n,t)$.\par

  We assume that $n$ is odd.
  Let $2k + 1$ be the greatest common divisor of $n$ and $t$. 
  In order to construct a Hamilton cycle in DP$(n, t)$, we define paths $P_i, Q_i, R_i, S_i$ for all $i \in \mathbb{Z}_{2k+1}$.
  \begin{align*}
    P_i &\colon u_{a_i+t} x_{a_i+t}  x_{a_i+t+1}  x_{a_i+t+2}  \dots  x_{a_{i+2}+t-1}  u_{a_{i+2}+t-1} \\
    Q_i &\colon v_{a_i}   y_{a_i}    y_{a_i+1}    y_{a_i+2}    \dots  y_{a_{i+2}-1}    v_{a_{i+2}-1} \\
    R_i &\colon u_{a_{i+1}+t-1}  v_{a_{i+1}+2t-1}  u_{a_{i+1}+3t-1}  \cdots v_{a_i} \\
    S_i &\colon v_{a_{i+1}-1}  u_{a_{i+1}-t-1}  v_{a_{i+1}-2t-1}  \cdots u_{a_i+t}
  \end{align*}
  where $a_0, a_1, a_2, \dots ,a_{2k} \in \mathbb{Z}_{2k+1}$ satisfy the following conditions
  \begin{align*}
    &\forall i \in \mathbb{Z}_{2k+1}, a_i \equiv i \ (\text{mod}\ 2k+1) \\
    &0 \leq a_0 < a_2 < a_4 < \dots < a_{2k} < a_1 < a_3 < a_5 < \dots < a_{2k-1} < n
  \end{align*}
  For instance, if $a_0=0, a_2=2, a_4=4, \dots, a_{2k}=2k, a_1=2k+2, a_3=2k+4, a_5=2k+6, \dots ,a_{2k-1}=4k$,
  the above conditions are met.
  Joining the paths in the following way gives a Hamilton cycle in DP$(n,t)$.
  \begin{align*}
    ((S_0 - P_0) - (R_1 - Q_1) - (S_2 - P_2) - (R_3 - Q_3) - \cdots \\
        \cdots - (R_{2k-1} - Q_{2k-1}) - (S_{2k} - P_{2k})) - \\
    - ((R_0 - Q_0) - (S_1 - P_1) - (R_2 - Q_2) - (S_3 - P_3) - \cdots \\
        \cdots - (S_{2k-1} - P_{2k-1}) - (R_{2k} - Q_{2k}))
  \end{align*}

\begin{figure}[!tpb]
    \begin{center}
      \includegraphics[height=0.44\textheight]{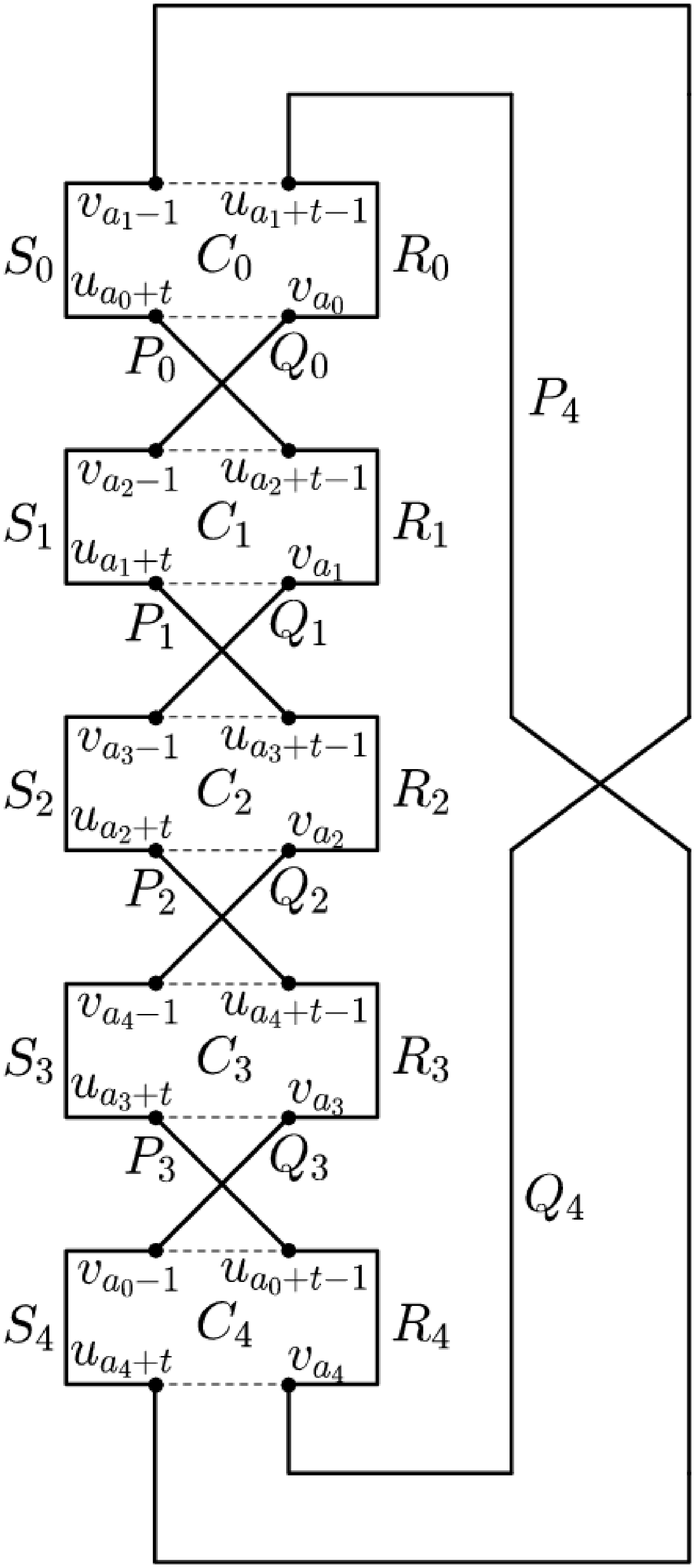}
      \caption{A Hamilton cycle in DP$(n, t)$ $(2k + 1 = 5)$}
      \label{fig:cycle}
    \end{center}
  \end{figure}

  An example of a Hamilton cycle in DP$(n, t)$ is shown in Figure \ref{fig:cycle}.

\section{Proof of Theorem \ref{thm}}
  In this section, we prove that the cycle described in the previous section contains all vertices of DP$(n,t)$ exactly once for all odd integers $n \geq 3$.
  Let $G$ be DP$(n, t)$ and $2k+1$ be the greatest common divisor of $n$ and $t$.
  \par
  Firstly, we prove that paths $Q_0, Q_1, \dots, Q_{2k}$ contain all of $V_y(G)$.
  \begin{align*}
    \bigcup_{i=0}^{2k} V_y(Q_i)
    = &\left(\bigcup_{i=0}^{k} V_y(Q_{2i}) \right) \cup \left(\bigcup_{i=0}^{k-1} V_y(Q_{2i+1}) \right) \\
    = &( \{y_{a_0}, y_{a_0+1}, \dots ,y_{a_2-1}\} \cup\\
      &\ \{y_{a_2}, y_{a_2+1}, \dots ,y_{a_4-1}\} \cup\\
      &\ \{y_{a_4}, y_{a_4+1}, \dots ,y_{a_6-1}\} \cup\\
      &\ \ \vdots\\
      &\ \{y_{a_{2k}}, y_{a_{2k}+1}, \dots ,y_{a_1-1}\}) \cup \\
      &( \{y_{a_1}, y_{a_1+1}, \dots ,y_{a_3-1}\} \cup\\
      &\ \{y_{a_3}, y_{a_3+1}, \dots ,y_{a_5-1}\} \cup\\
      &\ \{y_{a_5}, y_{a_5+1}, \dots ,y_{a_7-1}\} \cup\\
      &\ \ \vdots\\
      &\ \{y_{a_{2k-1}}, y_{a_{2k-1}+1}, \dots ,y_{a_0-1}\})\\
    = &\{y_m \mid m \in \mathbb{Z}_n\}\\
    = &V_y(G)
  \end{align*}
  \par
  According to the above equation and the definitions of $P_i$ and $Q_i$,
  we can prove that paths $P_0, P_1, \dots, P_{2k}$ contain all of $V_x(G)$.\par
  
  Secondly, we will prove that paths $R_0,R_1,\dots,R_{2k},S_0,S_1,\dots,S_{2k}$ contain all of $V_u(G) \cup V_v(G)$.
  We define cycles $C_i$ in DP$(n, t)$ for all $i \in \mathbb{Z}_{2k+1}$.
  \begin{align*}
    C_i \colon u_i v_{i+t}  u_{i+2t}  v_{i+3t}  \cdots  u_{i+(p-1)t} v_i u_{i+t}  v_{i+2t}  u_{i+3t}  \cdots  v_{i+(p-1)t}  u_i
  \end{align*}
  Note that odd integers $p$ and $q$ satisfy $n = p(2k + 1)$ and $t = q(2k + 1)$.
  For all $i \in \mathbb{Z}_{2k+1}$, $C_i$ consists of paths $D_i \colon u_i v_{i+t} u_{i+2t} v_{i+3t} \cdots u_{i+(p-1)t}$ and
  $E_i \colon v_i u_{i+t} v_{i+2t} u_{i+3t} \cdots v_{i+(p-1)t}$.
  Since $p$ is odd, the last vertex of $D_i$ is not $v_{i+(p-1)t}$ but $u_{i+(p-1)t}$.
  By symmetry, the last vertex of $E_i$ is $v_{i+(p-1)t}$.
  In addition, $u_{i+(p-1)t}$ and  $v_{i+(p-1)t}$ are respectively adjacent to $v_i$ and $u_i$ since $pt = pq(2k + 1)$ is a multiple of $n$.
  Observe that $u_i, u_{i+t}, u_{i+2t}, u_{i+3t}, \dots ,u_{i+(p-1)t}$
  contain no two same vertices since $pt$ is the least common multiple of $n$ and $t$.
  Hence $v_i, v_{i+t}, v_{i+2t}, v_{i+3t}, \dots ,v_{i+(p-1)t}$ also contain no two same vertices.
  \par

  We show that cycles $C_0, C_1, \dots, C_{2k}$ contain all of $V_u(G) \cup V_v(G)$.
  \begin{align*}
       \bigcup_{i=0}^{2k} V_u(C_i)
    &= \bigcup_{i=0}^{2k} \{u_{i+jt} \mid 0 \leq j < p\} \\
    &= \bigcup_{i=0}^{2k} \{u_{i+jq(2k+1)} \mid 0 \leq j < p\} \\
    &= \bigcup_{i=0}^{2k} \{u_{i+j(2k+1)} \mid 0 \leq j < p\} \\
    &= \bigcup_{i=0}^{2k} \{u_m \mid m \in \mathbb{Z}_n, m \equiv i\ (\text{mod}\ 2k+1)\} \\
    &= \{u_m \mid m \in \mathbb{Z}_n\}
  \end{align*}
  By symmetry, we have
  \begin{align*}
    \left(\bigcup_{i=0}^{2k} V_u(C_i)\right) \cup \left(\bigcup_{i=0}^{2k} V_v(C_i)\right)
    &= \{u_m \mid m \in \mathbb{Z}_n\} \cup \{v_m \mid m \in \mathbb{Z}_n\} \\
    &= V_u(G) \cup V_v(G)
  \end{align*}

  Observe that both $R_i$ and $S_i$ are subgraphs of $C_i$.
  For all $i \in \mathbb{Z}_{2k+1}$, $R_i$ and $Q_i$ share no vertex and contain all vertices in $C_i$
  since the first vertex of $R_i$ and the last vertex of $R_i$
  are respectively adjacent to the first vertex of $S_i$ and the last vertex of $S_i$.
  Therefore, paths $R_0,R_1,\dots,R_{2k},S_0,S_1,\dots,S_{2k}$ contain all of $V_u(G) \cup V_v(G)$.
  This completes the proof of Theorem \ref{thm}. \sq



\begin{thebibliography}{1}
\expandafter\ifx\csname url\endcsname\relax
  \def\url#1{\texttt{#1}}\fi
\expandafter\ifx\csname urlprefix\endcsname\relax\def\urlprefix{URL }\fi
\expandafter\ifx\csname href\endcsname\relax
  \def\href#1#2{#2} \def\path#1{#1}\fi

\bibitem{Wat}
M.~E. Watkins, A theorem on {T}ait colorings with an application to the
  generalized {P}etersen graphs, J. Combin. Theory 6 (1969) 152--164.

\bibitem{Cas}
F.~Castagna, G.~Pins, Every generalized {P}etersen graph has a {T}ait coloring,
  Pacific J. Math 40 (1972) 53--58.

\bibitem{Als}
B.~Alspach, The classification of {H}amiltonian generalized {P}etersen graphs,
  J. Combin. Theory Ser. B 34 (1983) 293--312.

\bibitem{Fu}
Y.~Y. Xueliang~Fu, B.~Jiang, On the domination number of generalized {P}etersen
  graphs ${P}(n,2)$, Discrete Math 309 (2009) 2445--2451.

\bibitem{Jin1}
J.-X. Zhou, Y.-Q. Feng, Cubic vertex-transitive non-{C}ayley graphs of order
  $8p$, Electron. J. Combin 19 (2012) {\#}P53.

\bibitem{Jin2}
J.-X. Zhou, Y.-Q. Feng, Cubic bi-{C}ayley graphs over abelian groups, European
  J. Combin 36 (2014) 679--693.

\bibitem{Kutnar}
K.~Kutnar, P.~Petecki, On automorphisms and structural properties of double
  generalized {P}etersen graphs, Discrete Math 339 (2016) 2861--2870.

\end{thebibliography}
\end{document}